\documentstyle{article}
\def \C{{\Bbb C}}
\def \ad{{\rm Ad}}
\def \calf{{\cal F}}
\def \P{{\Bbb P}}

\def \N{{\cal N}}
\def \S{{\Sigma}}
\def \R{{\Bbb R}}
\def \H{{\Bbb H}}

\def \min{{\setminus}}
\def \l{{\lambda}}
\def \Z{{\Bbb Z}}

\def \tag{{\tilde{\frak{g}}}}
\def \ag{{\frak{g}}}

\def \o{{\omega}}
\def \K{{K\"{a}hler\ }}
\def \L{{\cal L}}
\def \M{{\cal M}}
\def \CC{{\cal C}}
\def \fd{{\bullet}}

\def \MN{{{\cal M}_{\cal N}}}

\def \i{{\sqrt{-1}}}
\def \d{{\partial}}
\def \G{{\cal G}}
\def \g{{\gamma}}
\def\min{\setminus}

\def \proof{{\noindent{\it Proof.\ \ }}}
\newtheorem{Th}{THEOREM}[section]

\newtheorem{prop}[Th]{PROPOSITION}

\newtheorem{lem}[Th]{LEMMA}

\title{Moduli of flat bundles on open K\"{a}hler manifolds.}
\author{Jean-Luc Brylinski${}^1$,\ \ Philip A. Foth}
\begin{document}
\maketitle
\input amssym.def
\begin{abstract}
We consider the moduli space $\MN$ of flat unitary connections on an open
K\"{a}hler manifold $U$ (complement of a divisor with normal crossings) with
restrictions on their monodromy transformations.  Using intersection cohomology
with degenerating coefficients we construct a natural symplectic
form  $F$ on $\MN$. When $U$ is quasi-projective we prove that $F$ is actually
a \K form.  \end{abstract}

\tableofcontents
\footnotetext[1]{The first author was supported in part by NSF grant
DMS-9504522.}

\section{Introduction}
\setcounter{equation}{0}

Let $X$ be a compact K\"{a}hler manifold and $D$ a divisor on $X$ with normal
crossings. There exists a moduli space $\MN$ of flat irreducible unitary bundles
on $U=X\min D$ such that the monodromy transformation around each smooth
irreducible component of $D$ lies in a prescribed conjugacy class in $U(N)$.  We
develop a theory of deformations of representations of a discrete group relative
to our situation. We obtain a condition (similar to Goldman-Millson \cite{GM})
when $\MN$ has a manifold structure and further we work under this assumption.
If we pick a representation $\rho$ of $\pi_1(U)$ satisfying these conditions,
then we get a local system $\tag$ on $U$ associated with the representation
$\rho$.  Since $\tag$ has singularities on $X$ (along $D$) it is quite natural
to express the tangent space $T_{\rho}\MN$ in terms of the intersection
cohomology groups of $X$ with coefficients in $\tag$. Lemma \ref{lem:l32} shows
that this tangent space identifies with the group $IH^1(X, \tag)$. We introduce
a natural $2$-form on $\MN$ as a pairing $IH^1(X, \tag)\times IH^1(X,
\tag)\to\R$.

To see that our form is actually closed, we map our manifold $\MN$ to
an infinite-dimensional affine space such that its tangent space at any point
consists of $L_2$ forms on $X$ having certain additional properties. This affine
space admits a constant coefficient $2$-form such that its pull-back to $\MN$ is
given by a pairing in $L_2$ cohomology.  We further use the isomorphism between
intersection cohomology and $L_2$ cohomology constructed by
Cattani-Kaplan-Schmid \cite{CKS}, and Kashiwara-Kawai \cite{KK}.  As an
auxiliary tool we introduce the notion of $L_2$ vector bundle, which seems to be
of independent interest.

We prove \begin{Th} The moduli space $\MN$ is symplectic. \label{Th:t11}
\end{Th}

Intersection cohomology enjoys
such properties as Poincar\'{e} duality and the Hard Lefschetz theorem. This
follows from an unpublished work of Deligne and is partially explained by Zucker
in \cite{Zu}. The point is that the Hard Lefschetz theorem holds true
for the $L_2$ cohomology groups, because they are finite-dimensional, which
again is a consequence of the isomorphism between $L_2$  and intersection
cohomologies. The Hard Lefschetz theorem provides us with an isomorphism
$IH^{d-j}(X, \tag)\simeq IH^{d+j}(X, \tag)$, where $d=\dim_{\C}X$. This allows
us to see that our $2$-form is non-degenerate.

We also prove \begin{Th} When $X$ is projective, the moduli space $\MN$ is \K.
\label{Th:t12} \end{Th} When  $X$ is a curve we use the identification of $\MN$
with the moduli space of stable parabolic vector bundles given by Mehta-Seshadri
\cite{MS} to get a \K structure on the space $\MN$. Then we reduce the case of
general projective manifold $X$ to the case of a curve by taking appropriate
number of hyperplane sections of $X$.

The result was known in some special cases.  For example, in the compact case
(when $D$ is empty), the symplectic structure on the moduli space appears in the
works of Atiyah-Bott \cite{AB}, Goldman  \cite{G}, and Karshon \cite{Kar}.  In
the case of a Riemann surface with punctures the symplectic structure on the
moduli space was described by Atiyah \cite{A3}, Biquard \cite{BiqTh} (see also
\cite{Biq2}), and Witten \cite{W}. It was also the subject of a paper by
Biswas-Guruprasad \cite{BG}; a proof in this case using group cohomology is due
to Guruprasad-Huebschmann-Jeffrey-Weinstein \cite{GHJW}.

We would like to thank P. Deligne for his useful comments.

\section{Description of the moduli space}
\setcounter{equation}{0}

Let $X$ be a compact K\"ahler manifold of complex dimension $d$ endowed with a
K\"{a}hler form $\l$, and let $D$ be a divisor on $X$ with normal crossings such
that $D=\cup_{i=1}^rD_i$ is a decomposition of $D$ into the union of smooth
irreducible complex analytic subvarieties. Let $G$ be a compact connected Lie
group and $\ag$ its Lie algebra with fixed non-degenerate invariant bilinear
form $B\langle\ ,\ \rangle$.  From now on we fix the set $\N=(\CC_1, \CC_2, ...,
\CC_r)$ of $r$ conjugacy classes in $G$.  On the complement of $D$ we also have
a K\"{a}hler form - the restriction of the form $\l$ so that one has an open
K\"{a}hler manifold $U:=X\setminus D$. We notice that any quasi-projective
smooth algebraic variety can be obtained this way.

Take a base point $b\in U$ and denote by $\pi_1(U)=\pi_1(U,b)$ the fundamental
group of $U$ with base point $b$ and by $\pi_1(X)=\pi_1(X,b)$ the fundamental
group of $X$ with the same base point. Let us define ${\bar\M}_{\N}$ - the
moduli space of flat $G$-bundles on $U$ such that the monodromy transformation
around $D_i$ lies in $\CC_i$. Later on we will provide the reader with some
examples when ${\bar\M}_{\N}$ is actually smooth.

The moduli space $\M$ of flat $G$-bundles over $X$ identifies with the space
$Hom(\pi_1(X), G)/G$, where the group $G$ acts by conjugation.  As it is
well-known \cite{G} \cite{Kar} the space $\M$ admits a natural symplectic
structure.  Let $\rho$ be a smooth point on $\M$. Then the Zariski tangent space
at the class of $\rho$ is $T_{[\rho]}\M=H^1(\pi_1(X), \ag)$, where $\ag$ is
considered as a $\pi_1(X)$-module via the adjoint representation followed by
$\rho$.  Let $\tag$ stand for the local system associated with the
representation $\rho$.  Due to the isomorphism $H^1(\pi_1(X), \ag)\simeq H^1(X,
\tag)$ and the natural map $H^2(\pi_1(X), \ag)\to H^2(X, \tag)$ we get the
operation $$ \begin{array}{ccccc} H^1(\pi_1(X), \ag) & \cup & H^1(\pi_1(X), \ag)
& \to  & H^2(\pi_1(X), \ag) \\ \parallel & {} & \parallel & {} & \downarrow \\
H^1(X, \tag) & \cup & H^1(X, \tag) & \to & H^2(X, \tag) \end{array} $$ Then we
compose $$ B\langle\cdot,\cdot\cup\ [\l]^{d-1}\rangle:  H^1(X, \tag)\times H^1(X,
\tag) \to\R $$ to complete the pairing.

Our purpose is to construct a symplectic structure (and, moreover a K\"{a}hler
structure) on the moduli space ${\bar\M}_{\N}$ (at least on its smooth locus)
satisfying the following condition:  if $X$ is smooth and $\CC_i=1$ for each $i$
(under this assumption of course $\M={\bar\M}_{\N}$) then this form is the one
described above.  The moduli space ${\bar\M}_{\N}$ is the same as
$Hom_{\N}(\pi_1(U), G)/G$, where the subscript $\N$ means that  the prescribed
generator of $\pi_1(U)$ corresponding to a given loop $\g_i$ around $D_i$ goes
to $\CC_i\in G$. We will denote by $\MN$ the locus of ${\bar\M}_{\N}$
corresponding to the {\it irreducible} representations of $\pi_1(U)$.

\

\noindent{\it Remark.} The conjugacy classes $\CC_i$ may be chosen in such a
fashion that the moduli space is empty (on $\C\P^1$ with one puncture one can
take $\CC\ne Id$, see also \cite{SPM}).  Also we notice that the obvious map
$$\pi_1(U)\to\pi_1(X) $$ is surjective.

\

\begin{prop} The Zariski tangent space to the moduli space $\MN$ at the point
$\rho$ is $$ T_{[\rho]}\MN=Ker(H^1(\pi_1(U), \ag)\to\prod_i H^1(\Gamma_i, \ag)),
$$ where $\Gamma_i\simeq\Z$ is generated by the class of a loop encircling
$D_i$.  \label{prop:p21} \end{prop} \proof We shall repeatedly use the fact that
for every connected manifold $Z$ the two groups $H^1(Z, \tag)$ and
$H^1(\pi_1(Z), \ag)$ are canonically isomorphic.  To understand the tangent
space $T_{[\rho]}\MN$ we recall the well-known fact that it is a subspace of the
tangent space $T_{[\rho]}(Hom(\pi_1(U),G)/G)$.  Also we note that the tangent
space to the conjugacy class $\CC_g$ of $g\in G$ is the subspace of $\ag$ given
as the range of $Ad(g)-1$ if we identify as usual the tangent space $T_gG$ with
$\ag$ via the action of the left translation by $g$. Now if $\gamma_i$ is the
class of a loop encircling  $D_i$ and if $\Gamma_i\subset\pi_1(U)$ is the cyclic
subgroup generated by $\gamma_i$ then the cohomology group $H^1(\Gamma_i, \ag)$
is the cokernel of the map $Ad(g)-1:\ \ag\to\ag$, so the equality follows.
$\bigcirc$

\

The space $Ker(H^1(\pi_1(U), \ag)\to\prod_i H^1(\Gamma_i, \ag))$ is called
{\it the parabolic cohomology} of $X$ with coefficients in the local system
$\tag$ on $U$ associated to the action of $\pi_1(U)$ on $\ag$.

We will briefly discuss some of the known results for punctured Riemann
surfaces. Here $D_1, ...,D_r$ are just distinct points on $\S$. Let
$G_{\C}=SL(n, \C)$, $SO(n, \C)$, or $Sp(2n, \C)$ and let $G$ be its standard
maximal compact subgroup.  Let $\C^r$ be the standard representation space of
$G_{\C}$ (so that $r=n$ or $2n$) and consider $V=\oplus_{i=1}^r\wedge^i\C^r$
which is naturally a representation space of $G_{\C}$ too. We say that an
element $A\in G_{\C}$ (or its conjugacy class) satisfies {\it property P} if two
stabilizers have the same dimension:  $\dim(V^{A_s})=\dim(V^T)$, where $T$ is a
maximal torus in $G_{\C}$ (or in $G$) and $A=A_sA_n$ is a Jordan decomposition
of $A$ into the product of commuting unipotent and semisimple elements. In terms
of the eigenvalues $\l_1, ..., \l_n$ of $A$ property P means that the product
$\l_{i_1}\cdots\l_{i_k}\ne 1$ for any $i_1<\cdots <i_k$ and $k<n$.  For more
 details we refer to a paper of the second author \cite{F} where it was shown
that in the case of a Riemann surface with one puncture the moduli space of
flat $SL(n,\C)$- (or $SU(n)$-) bundles is smooth if and only if the monodromy
transformation around the puncture has property P. It was also shown that if
$G_{\C}=SO(n, \C)$, or $Sp(2n, \C)$, then the moduli space has at worst quotient
singularities (i.e. it is a quotient of a smooth manifold by an action of a
finite group).

\section{Deformations of representations of a discrete group}
\setcounter{equation}{0}

We have to make an additional assumption in order to guarantee that each
point $[\rho]\in\MN$ (corresponding to an irreducible representation $\rho$) is
a smooth point and the tangent space $T_{[\rho]}\MN$ to the {\it manifold} is
given by Proposition 2.1. For this it is enough to show that every
infinitesimal deformation (an element of the Zariski tangent space) is tangent
to an analytic path in $\MN$. We recall a theorem of M. Artin \cite{Art}, which
asserts that an infinitesimal deformation is tangent to an analytic path if and
only if there exists a formal power series deformation with the infinitesimal
deformation as its leading term. In the compact \K case, the criterion is
established by Goldman-Millson \cite{GM}. For example, it is enough to have
$H^2(\pi_1(X), \ag)=0$. Analogously to their methods we will show that
in our situation a sufficient condition is the vanishing of
the second relative cohomology group $H^2(\pi_1(U), (\Gamma_i), \ag)$.

Therefore, for the rest of the paper we will assume that we are dealing
with the situation when $H^2(\pi_1(U), (\Gamma_i), \ag)=0$ and thus we have no
obstruction for representing the actual tangent space as the kernel of the map
$H^1(\pi_1(U), \ag)\to\prod_iH^1(\Gamma_i, \ag)$. In all other situation we only
get "formal symplectic structure" and it is a separate problem which will be
treated elsewhere to understand its actual meaning.

In the present section we will establish all the results in the following
generality.  Let $\pi$ be a discrete finitely generated group and let
$(\Gamma_i)_{i=1}^r$ be a system of its subgroups such that $\Gamma_i\simeq\Z$
is generated by an element $\g_i$. We consider the representation variety
$\MN=Hom(\pi, G)_{\N}/G$ ($G=U(N)$), consisting of classes of such group
homomorphisms that $Im(\pi)$ is not contained in any proper parabolic subgroup
of $G$ (irreducibility condition) and $Im(\g_i)\in\CC_i$. We saw that the
Zariski tangent space to a point $[\rho_0]$ is given by
$$T_{[\rho_0]}\MN=Ker[H^1(\pi, \ag)\to\oplus_iH^1(\Gamma_i, \ag)].$$ Our task is
to prove \begin{prop} If $H^2(\pi, (\Gamma_i), \ag)=0$ then for any $\eta\in
T_{[\rho_0]}\MN$ there exists a formal power series $\rho_t$ of $\rho_0$
representing $\eta$.  \label{prop:pg} \end{prop} We devote the rest of the
section to the proof of this statement.

Let us have $$\rho_t(g)=\rho_0(g)\exp(\sum_{i=1}^{\infty}f_i(g)t^i),$$ where
$f_i: \pi\to\ag$ is a group $1$-cochain and $t$ is a formal parameter. We recall
that $\ag$ is a $\pi$ - module via the adjoint representation followed by
$\rho_0$. We have the following two conditions to satisfy. First, $\rho_t$ is a
group homomorphism, therefore $\rho_t(g_1g_2)=\rho_t(g_1)\rho_t(g_2)$. Secondly,
the condition of mapping $\g_i$ to fixed conjugacy classes is
written as \begin{equation}\rho_t(\g_j)=\exp(\sum_{i=1}^{\infty}c_i^jt^i)
\rho_0(\g_j) \exp(-\sum_{i=1}^{\infty}c_i^jt^i),\label{eq:e30}\end{equation}
where $c_i^j\in\ag$. Next we make a change of notation $h_i(g)=f_i(g^{-1})$ and
using the Campbell-Hausdorff formula spell out these conditions. To the first
order of $t$ we have $$\d h_1 (g_1,
g_2)=\ad\rho_0(g_1)h_1(g_2)-h_1(g_1g_2)+h_1(g_1)=0, $$
$$h_1(\g_j)=\ad\rho_0(\g_j)c_1^j-c_1^j=\d c_1^j(\g_j),$$ where $c_1^j$ is
considered as a zero-cochain. As one sees, these equations are equivalent to the
fact that $h_1\in Ker[H^1(\pi, \ag)\to\oplus_iH^1(\Gamma_i, \ag)]$.

Our purpose is to find such $h_2, h_3, ...$ and $c_2^j, c_3^j, ...$ which
satisfy those two conditions. We will do this by the induction process.
Let us first show explicitly the existence of such $h_2$ and $c_2^j$.
To the second order of $t$ we have \begin{equation}
\begin{array}{c} \d h_2(g_1,g_2)=-{1\over
2}[\ad\rho_0(g_1)h_1(g_2),h_1(g_1)], \\ h_2(\g_j)-\d c_2^j=-{1\over
2}[\ad\rho_0(\g_j)c_1^j, c_1^j].
\end{array} \label{eq:e1} \end{equation}

We recall that given a group $\pi$ and a system of its subgroups $\Gamma_j$
together with the restriction maps $Map(\pi, \ag)\to Map(\Gamma_j, \ag)$
the relative cochain complex is defined as the cone of the system of
maps of complexes $R_j: C^{\fd}(\pi, \ag)\to C^{\fd}(\Gamma_j, \ag)$. By
definition, $$Cone_{\ag}^k(R_1, ..., R_r)=C^k(\pi, \ag)\oplus
\bigoplus_jC^{k-1}(\Gamma_j, \ag)$$ with the differential $(-\d, R_j + \d)$.

Therefore, we would be able to find such $h_2$ and $c_2^j$ if the relative
group $H^2(\pi, (\G_i), \ag)$ vanishes and the right hand side of (\ref{eq:e1})
is a relative $2$-cocycle. One easily checks that the cocycle condition is
satisfied: $$\d[\ad\rho_0(x)h_1(y), h_1(x)]\ (g_1, g_2, g_3)=0, \ \ g_1, g_2,
g_3\in\pi, \ \ {\rm and}$$ $$[\ad\rho_0(\tau_1)h_1(\tau_2),
h_1(\tau_1)]=\d[\ad\rho_0(x)c_1^j, c_1^j]\ (\tau_1, \tau_2), \ \ \tau_1,
\tau_2\in\Gamma_j.$$ (To verify the equalities we use that $\d h_1(g_1, g_2)=0$
and $h_1(\g_j)=\d c_1^j(\g_j)$.)

This procedure serves several purposes. First, we get an idea that the
obstruction for an element of Zariski tangent space to be tangent to an analytic
path lies in the second relative group cohomology. Besides, we notice that in
all successive steps we will deal with the system of the form \begin{equation}
\begin{array}{c}\d h_{k+1} (g_1, g_2, g_3)=F(h_1, ..., h_k) \\
h_{k+1}(\g_j)-\d c_{k+1}^j(\g_j)=H(c_1^j, ..., c_k^j) \end{array} \label{eq:e2}
\end{equation} Also, we see that in order to show that the right hand side of
these equations is a relative cocycle, we can restrict ourselves to the case of
just one subgroup $\Gamma\subset\pi$ generated by an element $\g$.

We need the following simple  \begin{lem} Let $\pi$ be a finitely generated
discrete group and let $\g\in\pi$. There exists a free group $\calf $ with
$g\in{\cal F}$  and a surjective homomorphism $\phi: \calf\to\pi$ such that
$\phi(g)=\g$ and $g$ is an element of a basis of $\calf $. \end{lem} \proof Let
$\calf'$ be any free group such that there is a surjective map $\phi':
\calf'\to\pi$. If $(g_2, ..., g_l)$ is a basis for $\calf'$ then we consider a
free group $\calf $ obtained from $\calf'$ by adding one generator $g$. Then
$(g, g_2, ..., g_l)$ is a basis of $\calf $ and we let $\phi(g)=\g$ and
$\phi(g_i)=\phi'(g_i)$ for $2\le i\le l$. $\bigcirc$

\

Now we make the inductive step. Let $\calf$ be such a free group that satisfies
the condition of the above lemma for our group $\pi$ and its subgroup $\Gamma$
generated by $\g$. Let $\Gamma'\simeq\Z$ be the subgroup of $\calf$ generated
by the element $g$ from the above lemma. (We notice that $H^2(\calf, \Gamma',
\ag)=0$.) Let us lift the system of equations (\ref{eq:e2}) to the free group
$\calf$. We were able to find a solution of this system up to order $k$.
The character variety of a free group in $l$ generators $(g, g_2, ..., g_l)$
such that the image of $g$ goes to $\CC\subset G$ is just
$$\CC\times\underbrace{G\times G\times \cdots \times G}_{l-1}$$ and it is
non-singular.

Thus there is no obstruction to finding $\rho_t:\calf\to G$ up to order $k+1$,
inducing the given $k$-th order formal homomorphism, as well as $C_{k+1}^j$ such
that \ref{eq:e30} helds up to $T^{k+1}$. We conclude that in the free group
$\calf$ it is possible to find a solution of the lift of (\ref{eq:e2}). This
means that the right hand side not only of the lift of (\ref{eq:e2}) but also of
(\ref{eq:e2}) itself is a relative cocycle.  It is a relative coboundary as
well, because we assumed the vanishing of the group $H^2(\pi, (\Gamma_i), \ag)$.
Therefore there exist $h_{k+1}$ and $c_{k+1}^j$ satisfying the equations
(\ref{eq:e2}), which completes the inductive step. This finishes the proof of
Proposition \ref{prop:pg}.

We refer to results \cite{KM} of Kapovich-Millson for another approach to the
relative deformation theory, where they work with differential graded algebras
of differential forms on a manifold.

\section{Intersection cohomology and the construction of the $2$-form}
\setcounter{equation}{0}

We will construct a non-degenerate $2$-form on the space $\MN$.  We first
 consider the dimension $1$ case.  The case of a Riemann surface $\S$ with
punctures is quite simple because we know the explicit structure of the
fundamental group and   $\S$ is a $K(\pi, 1)$.  The latter property allows us to
 conclude that $H^i(U, \tag)=H^i(\pi_1(U), \ag)$.  Consider the exact sequence
$$ \cdots\to H^1(Cone,\ag)\to H^1(\pi_1(U), \ag)\to \oplus_iH^1(\Gamma_i,
\ag)\to H^2(Cone_{\ag})\to 0, $$ where $Cone_{\ag}$ is the mapping cone for the
morphism of complexes $$ C^{\fd}(\pi_1(U), \ag)\to \bigoplus_iC^{\fd} (\Gamma_i,
\ag).  $$ Similarly one can define $Cone_{\Z}$, $Cone_{\R}$, etc.  If we apply
the bilinear form $B\langle ,\rangle$ together with the pairing in cohomology
then we get a map $$ H^i(Cone_{\ag})\times H^j(Cone_{\ag})\to
H^{i+j}(Cone_{\R}).  $$ \begin{lem} $H^2(Cone_{\Z})\simeq\Z$.  \end{lem} \proof
Let $\Delta_i$ be a small disk centered at $i$-th marked point.  Then
$\Delta_i^*$ be obtained from $\Delta_i$ by removing the marked point.  We have
the class $\gamma_i'\in\H^1 (\Gamma_i,\Z)\simeq H^1(\Delta_i^*, \Z)$, which can
be defined as follows. In terms of group cohomology it corresponds to the map
$\phi: \Gamma_i\to\Z$ sending $\gamma_i$ to $1$. Under the above isomorphism
this class goes to the class defined in $H^1(\Delta_i^*, \Z)$ by the loop
$\gamma_i$.

Our point is that $$ H^2(Cone_{\Z})=Coker(H^1(\pi_1(U), \Z)
\stackrel{\alpha}{\to}\oplus_iH^1(\Gamma_i, \Z)); $$ the cokernel of the map
$\alpha$ is $\Z$ and is generated by the class of the element $$ (\g_1', ...,
\g_r')\in\oplus_iH^1(\Delta_i^*, \Z) \simeq\oplus_iH^1(\Gamma_i, \Z). \ \bigcirc
$$ This lemma gives us the idea how to construct a symplectic form on $\MN$:
$T_{[\rho]}\MN\times T_{[\rho]}\MN\to\R$, because $$
T_{[\rho]}\MN=Ker(H^1(\pi_1(U), \ag)\to\prod_iH^1(\Gamma_i, \ag))= $$ $$
Im(H^1(Cone_{\ag})\to H^1(\pi_1(U), \ag)) $$ and we have a natural pairing $$
H^1(Cone_{\ag})\times H^1(Cone_{\ag}) \stackrel{B\langle,
\rangle}{\longrightarrow}H^2(Cone_{\R})\to\R.  $$ This pairing is clearly
non-degenerate.

Now we move on to arbitrary dimensions.  Unfortunately, the usual cohomology
does not allow us to construct a non-degenerate pairing on the tangent space of
$\MN$ when the local system $\tag$ cannot be extended from $U$ to $X$, i.e. when
at least one of the conjugacy classes $\CC_i$ is not the class of the identity
of $G$.  That is the main reason why we need to use intersection cohomology
instead.

We introduce the following (canonical) filtration on $X$:  $X_0\subset
X_1\subset\cdots\subset X_{n-1}\subset X_n=X$, where $X_j$ consists of all
points that belong to at least $n-j$ smooth irreducible components $D_i$ of
$D$.  We notice that $X\setminus X_{n-1}=U$.  We will always work with
intersection cohomology for the middle perversity.

We  define the truncated complex $\tau_jC^{\fd}$ for a complex of sheaves
$C^{\fd}$. It is the complex which in degree $i$ is $$ C^i \ \ \ \ \ \ \ {\rm
if} \ \ i < j, $$ $$ Ker(C^j\to C^{j+1})\ \ {\rm if} \ \ i=j, $$ $$ 0   \ \ \
\ \ \ \ {\rm if} \ \ i>j.  $$ \begin{lem}
$Ker(H^1(\pi_1(U), \ag)\to\prod_iH^1(\Gamma_i, \ag))=IH^1(X, \tag)$.
\label{lem:l32} \end{lem} \proof An important thing is that the
codimension of $D$ in $X$ is $1$.  That is why in the complex
$i_!i^!IC^{\fd}_X(\tag)$ has its cohomology sheaves only in degree $\ge 2$. If
$j: U\hookrightarrow X$ and $i: D\hookrightarrow X$ are the inclusions then
one has an exact (distinguished) triangle of complexes (see \cite{B}, p.109)
\begin{equation} i_!i^!IC^{\fd}_X(\tag)\to IC^{\fd}_X(\tag)\to
Rj_*j^*IC^{\fd}_X(\tag)\to i_!i^!IC^{\fd}_X(\tag)[1],
\label{eq:e31} \end{equation} which gives rise to the long exact sequence in
cohomology (we need only a small part of it):  $$ \cdots\to H^1(X,
i_!i^!IC^{\fd}_X(\tag))\to H^1(X, IC^{\fd}_X(\tag))\to $$ $$ H^1(U,
IC^{\fd}_U(\tag))\to H^2(X, i_!i^!IC^{\fd}_X(\tag))\to
H^2(X,IC^{\fd}_X(\tag))\to\cdots $$ The vanishing mentioned above implies that
$H^1(X, i_!i^!IC^{\fd}_X(\tag))=0$.  Besides, from definition it follows that
$H^1(X, IC^{\fd}_X(\tag))=IH^1(X, \tag)$, and since $U$ is non-singular and the
local system $\tag$ is defined on $U$ one has $H^1(U,
IC^{\fd}_U(\tag))=H^1(U,\tag)$.  All this means that $$ IH^1(X, \tag)=Ker(H^1(U,
\tag)\to H^2(X, i_!i^!IC^{\fd}_X(\tag))).  $$ We will show that there is a
natural injection $$ H^2(X,
i_!i^!IC^{\fd}_X(\tag))\hookrightarrow\prod_iH^1(\Gamma_i, \ag), $$ and in view
of Proposition \ref{prop:p21} it is enough to prove the lemma.

First, we consider the case $r=1$, when the divisor $D$ is smooth and
irreducible.  We will use Deligne's construction of the intersection sheaf
complex as in \cite{GMII}. Let us consider a point $x\in D\subset X$ and a
small neighbourhood $V$ of $x$. The group $\Gamma=\pi_1(V\min(V\cap D))$ is
isomorphic to $\Z$ and is generated by a loop $\gamma$ encircling $D$. If $k:
V\hookrightarrow X$ is the inclusion of this neighbourhood, then the sheaf
$k_*\tag$ has a stalk $(k_*\tag)_x=\underbar H^0(\pi_1(V\min(V\cap D), G)$.
Passing to the derived functor and taking the sheaf cohomology we see that a
stalk $\underbar H^1(Rk_*\tag)_x=H^1(\Gamma, \ag)$, which is the cokernel of
$\rho(\gamma)-Id$ in $\ag$. (We recall that $\ag$ is a $\pi_1(U)$-module via
the adjoint representation followed by $\rho: \pi_1(U)\to G$.) For $i>1$ the
sheaf cohomology groups vanish:  $\underbar H^i(Rk_*\tag)=0$. We also notice
that $H^1(V\min(V\cap D), \tag)=H^1(\pi_1(V\min(V\cap D)), \ag)$.

Let ${\cal S}$ stand for a local system on $D$ with a stalk ${\cal
S}_x=H^1(\Gamma, \ag)$.  Now we have $IC^{\fd}_X(\tag)=\tau_0Rk_*\tag$, since
the canonical filtration simply amounts to $D\subset X$.  Thus we get an exact
triangle $$ \cdots\to IC^{\fd}_X\to Rk_*\tag\to i_*{\cal S}[-1]\to
IC^{\fd}_X[1]\to\cdots, $$ because the mapping cone for $IC^{\fd}_X\to Rk_*\tag$
is homotopic to the complex having in degree $1$ the sheaf $i_*{\cal S}$ and
nothing else.  Comparing this triangle with (\ref{eq:e31}) we see that
$i_!IC^{\fd}_X(\tag)\simeq {\cal S}[-2]$ and it means that $$ H^2(D,
i^!IC^{\fd}_X(\tag))=H^0(D, {\cal S}), $$ and it injects into $H^1(\Gamma, \ag)$
by definition of ${\cal S}$.

Similarly we can consider the case $r>1$. Here the canonical filtration is
$\cdots\subset M\subset D\subset X$, where the subvariety $M$ consists of all
points that belong to at least $2$ irreducible components of $D$.  For any
$x\in M$ we take its small enough neighbourhood $V$; the fundamental group
$\pi_1(V\min(D\cap V))\simeq\Z^m$, where $m$ is the number of intersecting
components of $D$ at $x$. If $l:M\hookrightarrow X$ is the inclusion of $M$
into $X$ then the complex $l_!l^!IC^{\fd}_X(\tag)$ has its cohomology sheaves
only in degree $\ge 3$, since the codimension of $M$ in $X$ is at least $2$.
Thus from the exact sequence $$ \cdots\to H^2(X, l_!l^!IC^{\fd}_X(\tag))\to
IH^2(X,\tag)\to IH^2(X\min M,\tag)\to\cdots $$ one concludes that we have an
injection $$ IH^2(X,\tag)\hookrightarrow IH^2(X\min M,\tag).  $$ Now the same
arguments as above show that the group $H^2(X, i_!i^!IC^{\fd}_X(\tag))$
naturally injects into the product $\prod_iH^1(\Gamma_i, \ag)$. $\bigcirc$

\

The above lemma proves that $T_{\rho}\MN\simeq IH^1(X, \tag)$, and
now we are ready to construct a $2$-form $F$ on the space $\MN$. The
idea is as follows. We have the class $[\l]\in H^2(X, \R)$ of the K\"{a}hler
form (when $X$ is a projective variety the class $[\l]$ is just the class of a
hyperplane section). Now the form $F$ is given by the pairing \begin{equation}
IH^1(X, \tag)\times IH^1(X, \tag)\to\R:\ \ \ \langle x,y\rangle=B\langle x,
y\cup\l^{d-1}\rangle. \label{eq:e32} \end{equation} Here we used the
intersection pairing described in \cite{GMII}.

Now we shall make a choice of Riemannian metric on $U$ with singularities
along $D$. Locally the description of our metric is as follows. Let $\Delta$
stand for the standard open disc in $\C$ given by $|z|<\varepsilon$ for some
positive number $\varepsilon$ and let $\Delta^*=\Delta\min 0$. The
intersection of a small neighbourhood of $x\in X$ with $U$ looks like
$(\Delta^*)^r\times\Delta^{d-r}$ if $r$ components of $D$ meet in $x$. In
local coordinates $z_1, ..., z_d$ these components are defined by equations
$z_1=0$, $z_2=0$, ..., $z_r=0$. On $\Delta$ we take the metric $dzd\bar{z}$ and
on $\Delta^*$ there is the Poincar\'{e} metric given in polar coordinates $r,
\theta$ by $$ {dr^2+(rd\theta)^2\over (r\ln r)^2}.  $$ So we assume that on $U$
we have a metric that is quasi-isometric to this one over any such open set
$(\Delta^*)^r\times\Delta^{d-r}$.

It is important to notice that the local system $\tag$ is unitary, since the
representation $\ag$ of $G$ is unitary. By a well-known theorem (\cite{CKS},
\cite{KK}) one has the isomorphism $IH^i(X, \tag)\simeq H^i_{(2)}(U,\tag)$
between the intersection cohomology and the $L_2$-cohomology with coefficients
in $\tag$. Therefore when $i=1$ this space carries a pure Hodge structure of
weight $1$, with only Hodge types $(0,1)$ and $(1, 0)$. This induces a
complex structure on the tangent space to $\MN$.

We refer the reader to \cite{BZ} for a survey of Hodge theory and its relation
with $L_2$  cohomology.

\section{The universal bundle}
\setcounter{equation}{0}
From now on we shall concentrate on the case when $G=U(N)$ and $G_{\C}=GL(N,
\C)$.  Let us consider a principal flat $G$-bundle $P$ over $U$; then there
exists canonically a holomorphic vector bundle $\bar V$ over $X$ such that $$
\bar{V}_{|U}=V:=P\times_{G}\C^N, $$ where $\C^N$ is the standard $G$-module.
The holomorphic bundle $\bar V$ is called the {\it Deligne extension} of $V$
(see \cite{De1}).  It has the property that the corresponding connection
$\nabla_0$ has at worst logarithmic singularities along $D$. The canonical
Deligne extension relies on fixing a unit interval: it is assumed that if
$\mu$ is an eigenvalue of $Res_{D_i}(\nabla_0)$ then $0\le Re(\mu)< 1$.

Over $U$ we have the Lie algebra bundle $End(V)$ and it extends to the
holomorphic bundle $End({\bar V})$ over $X$.

Let $Z$ be the subvariety of $Hom(\pi_1(U), G)$ consisting of all irreducible
representations with the monodromy transformation around $D_i$ lying in the
fixed conjugacy class $\CC_i$, so that $\MN=Z/G$.  Let us consider the vector
bundle ${\Bbb U}$ over $Z\times U$ given by $${\Bbb
U}=(\C^N\times\tilde{U}\times Z)/\pi_1(U), $$ where $\tilde{U}$ is
the universal covering of $U$ and the action of an element $a\in\pi_1(U)$ is
given by $$ a(x, \tilde{y}, \rho)=(\rho(a)(x), a(\tilde{y}), \rho), \ \
x\in\C^N, \ {\tilde y}\in \tilde{U}.  $$ Also one can construct the universal
bundle $E$ over the product $\MN\times U$:  $$ E=(\C^N\times{\tilde U}\times
Z)/\pi_1(U)\times G $$ $$ (a,g)(x, {\tilde y}, \rho)= (g\rho(a)g^{-1}(x),
a({\tilde y}), g\rho g^{-1}),\ \ \ g\in G.  $$ We will need the following
\begin{lem} For any $x\in D$ there exists a neighbourhood $V$ of $x$ such that
as $\rho$ varies in a connected component of $Z$  the local monodromy
representation $\Z^k\simeq\pi_1(V\min(V\cap D))\to G$ does not change (up to
conjugacy).  \end{lem} \proof Clearly one only has to show that the map $\xi$:
$$\begin{array}{c} \{\ k\ -{\rm\ tuples\ of\ commuting\ elements\ of\ } G\}
 /{\rm conjugacy} \\ \xi\downarrow \\ (G {\rm\ mod\ conjugacy})^k \end{array}$$
is a finite map. It is also important to notice that the source of this map is
Hausdorff. Let us exhibit this for $k=2$, because for $k>2$ the arguments are
the same.

Assume that for $a,b,c,d\in G$ such that $[a,b]=[c,d]=1$ one has
$\xi(a,b)=\xi(c,d)$. Due to the fact that we consider pairs up to conjugacy we
may assume that $a=c$. Moreover, $b$ is conjugate to $d$. We identify pairs
$(a,b)$ and $(a,d)$ if $b$ is conjugate to $d$ by means of an element from
$Z(a)$ - the centralizer of $a$. Next, we observe that each $G$-conjugacy
class intersects $Z(a)$ by only finitely many $Z(a)$-conjugacy classes.  In
fact, the number of these classes is bounded from above by the cardinality of
the Weyl group of $G$. $\bigcirc$

\

This lemma allows us to prove the next important result.  \begin{prop} There
exists a vector bundle $\tilde {\Bbb U}$ over $X\times Z$ extending ${\Bbb
U}\to U\times Z$ such that for every smooth point  $\rho\in Z$ the restriction
$\tilde{{\Bbb U}}_{|X\times\rho}$ is the holomorphic Deligne extension.
\end{prop} \proof The above lemma allows one to construct the bundle
$\tilde{{\Bbb U}}$ locally near $D\subset X$, because if one picks a point
$x\in D$ then in a neighbourhood $V\ni x$ the local representation of
$\pi_1(V\min(V\cap D))$ does not vary. So one can take the bundle over
$V\times Z$ as a pullback of the bundle over $V$ under the projection onto the
first coordinate. Now one notices that the Deligne construction is compatible
with restrictions to smaller open sets. It means that for neighbourhoods $V_1$
and $V_2$ of two points $x_1, x_2\in D$ respectively the restrictions of
${\tilde {\Bbb U}}$ onto $V_1$ and $V_2$ agree on the intersection $V_1\cap
V_2$. $\bigcirc$

\

Repeating the above arguments verbatim one can prove the existence of the
extension $\tilde{E}\to\MN\times X$ of the universal bundle $E$. In fact, the
bundle ${\tilde E}$ is the one we use in the rest of the paper. The bundle
$\tilde E$ can also be obtained using the quotient of ${\tilde{\Bbb U}}$ by the
action of $G$. We also notice that if the monodromy representation of
$\pi_1(U)$ is irreducible then the corresponding logarithmic connection is
stable in the sense of \cite{Nits}.

\section{Gauge group and $L_2$ bundles}
\setcounter{equation}{0}

In this section we introduce an "$L_2$ gauge group" $\G$, later we will use it
to prove Theorem \ref{Th:t11}.

As before we consider a (holomorphic) vector bundle $\bar V$ on $X$ which is the
Deligne extension of the vector bundle $V$ on $U$ with a fixed flat connection
$\nabla_0$.  The connection $\nabla_0$ has logarithmic poles along the divisor
$D$ and it corresponds to an irreducible unitary representation $\rho_0$ of
$\pi_1(U)$. This representation sends monodromy transformations around
irreducible components of $D$ to the prescribed set ${\cal N}$ of conjugacy
classes in $U(N)$.

Now we are ready to introduce the gauge group $\G$ that acts on the space of
connections on $\bar V$.  It is the set of smooth unitary automorphisms of the
bundle $V$ over $U$ such that if $g\in\G$ then $g^{-1}\nabla_0g$ is a unitary
$L_2$ $1$-form with coefficients in $End(V)$.

The topology on $\G$ is the coarsest topology satisfying the following three
conditions. First, for any compact subset $K$ of $U$, the map from $\G$ to the
Fr\'echet space of $C^{\infty}$ sections of $End(V)$ over $K$ should be
continuous.
Second, our topology is such that the distance function $G\times G\to\R$ defined
by  $${\rm dist}(g_1,g_2)=\sup_{x\in U, v\in V_x,
||v||=1}||g_1^{-1}g_2(v)-v||,$$ is continuous. The above supremum is
well-defined since every $g\in\G$ is bounded.  Finally, we require that in our
topology the function $\G\to\R$ given by the $L_2$-norm of the covariant
derivative of $g\in G$ is continuous.

To see explicitly the manifold structure, we will give a description of the
algebra $Lie(\G)$ and construct a map from a small neighbourhood of zero in
$Lie(\G)$ to a small neighbourhood of $1\in\G$.  The algebra $Lie(\G)$ consists
of bounded sections $u$ of $u(V)$ - the unitary Lie algebra bundle corresponding
to $V$ such that  $\nabla_0u$ is $L_2$. Now the map $Lie(\G)\to\G$ in the
neighbourhood of $0$ is just the Cayley transform $(\i I-A)(\i I+A)^{-1}$.

Let $\{ Y_i\}$ be a finite open cover of $D$ in $X$ by contractible open
polycylinders $\Delta^d$ so that $Y_i\simeq (Y_i\cap D)\times \Delta$. Next we
need to define {\it a hermitian $L_2$ vector bundle} over $X$. It is a triple
$(V, h, \{ B_i\})$, where $V$ is a vector bundle over $U$ with hermitian metric
$h$, and $B_i$ is a class modulo $\G$ of frames of $V$ over
$Y_i\cap U$ such that $B_i$ and $B_j$ agree over the intersection $Y_i\cap
Y_j\cap U$. If we have a refinement $\{ \tilde{Y}_j\}$ of $\{ Y_i\}$ with
similar properties then we naturally obtain a triple $(V, h, \{ \tilde{B_j}\})$,
defining the same hermitian $L_2$ vector bundle. Two hermitian $L_2$ vector
bundles are isomorphic if they are isomorphic over a common refinement. We call
elements of $B_i$ {\it local $L_2$-frames}.  We define a section of a hermitian
$L_2$ bundle to be a section of $V$ which is square integrable near $D$.

Let $x\in D$ be an arbitrary point of the divisor $D$ and let $Y$ be a small
open polycylinder with Poincar\'{e} metric containing $x$. Let $f$ be a smooth
compactly-supported function on $Y$. Then it is easy to see that $||df||$ in
Poincar\'{e} metric is bounded from above. This allows us to use partitions
of unity in the $L_2$-context.

Let $Z$ be a smooth manifold (possibly with boundary). We introduce the $L_2$
gauge group $\G'$ , which consists of those smooth maps $g\in Map(Z\times(X\min
D), G)$, which satisfy $$\sup_{z\in C}||\nabla_Xg(\cdot,z)||<\infty,$$ where $C$
is a compact subset of $Z$, $\nabla_X$ is the covariant derivative in
$X$-direction, and $||\cdot ||$ is the global $L_2$-norm on $X$. Let $p_2$ be
the projection $Z\times X\to X$. Then $p_2^{-1}(D)$ has real codimension $2$ in
$Z\times X$ (and still has normal crossings), and the notion of hermitian $L_2$
bundle over $Z\times X$ makes perfect sense. It is a triple $(E,h,\{
p_2^{-1}B_i\})$, where $E$ is a vector bundle over $Z\times(X\min D)$ equipped
with a hermitian metric $h$.   One can think of a hermitian $L_2$ bundle over
$Z\times X$ as of a family of hermitian $L_2$ bundles over $X$ varying smoothly
with $z\in Z$. The important step in our construction is the following Lemma
which has a well-known prototype in topological K-theory.

\begin{lem} Let $Z_1, Z_2$ be smooth contractible manifolds and let $\Phi$ be a
hermitian $L_2$ bundle over $Z_2\times X$. Let $a_0$ and $a_1$ be two smooth
homotopic maps $Z_1\to Z_2$. Then $(a_0\times Id)^*\Phi$ and $(a_1\times
Id)^*\Phi$ are isomorphic as hermitian $L_2$ vector bundles over $Z_1\times X$.
\end{lem}

\proof Let $I$ denote the unit interval and let $\L$ be a hermitian $L_2$ bundle
over $X\times I\times Z$. Let $K$ be a closed subset of $X\times I$ (actually,
we use $K=X\times\{ t\}$). The restriction  $\L_{K\times Z}$ is a hermitian
$L_2$ bundle over $K\times Z$. Locally, a section $s$ of this bundle is a vector-valued function on $K\times Z$ which by definition satisfies the following
estimates on its $L_2$-norms:  $$||s||<\infty, \ \ ||\nabla_X s||<\infty.$$
Therefore the Tietze extension theorem can be applied and for each $x\in X\times
I\times Z$ we may find an open set $Y\subset X\times I$ satisfying $x\in Y\times
Z$ and a section $s'$ of $\L$ such that $s$ and $s'$ agree over $(Y\cap K)\times
Z$. Since $X\times I$ is compact, we can find a finite cover by such open sets
and apply $L_2$ partition of unity to see that the section $s$ can be extended
to $X\times I\times Z$. To finish the proof, we apply the classical
argument as given on page 17 of \cite{AK}.  $\bigcirc$

Now we pick a small contractible neighbourhood $W\subset\MN$ of $\nabla_0$.  Let
$a_0, a_1: W\times X\to W\times X$ be two homotopic maps, where the map
$a_0$ is just the identity map and $a_1$ is the map $W\times X\to
\{\nabla_0\}\times X$ which induces identity on the second coordinate. The
bundle we are interested in is $\tilde E$ from Section 5. We see that the bundle
$a_1^*({\tilde E})$ is naturally isomorphic to the bundle $p_2^*{\bar V}$ over
$W\times X$, where $p_2:  W\times X\to X$ be the projection onto the second
coordinate. Now we apply the above Lemma to the maps $a_0$ and $a_1$ to get an
isomorphism $$\psi: {\tilde E}\simeq p_2^*{\bar V}$$ of hermitian $L_2$ bundles
over $W\times X$. For any $\rho\in W$ we let $\psi_{\rho}$ stand for the
restriction of this isomorphism to $\rho\times X$.

\section{The $2$-form is symplectic.}
\setcounter{equation}{0}

Our next goal is to show that the form $F$ constructed earlier is actually
closed. For this we shall construct another $2$-form $H$ on a bigger space such
that its closeness is apparent and then we shall see that our $2$-form $F$ is
equal to a pullback of the form $H$.

Let $Af$ be the affine space of all unitary connections $\nabla$ on $\bar V$
such that the difference $\nabla-\nabla_0$ is an $L_2$ $1$-form with
coefficients in $End({\bar V})$.

For each $\rho\in W$ we have a pull-back $\nabla_{\rho}=\psi_{\rho}^*\nabla_0$
of the connection $\nabla_0$ in $V$. The corresponding map $\beta: W\to Af$,
$\beta(\rho) = \nabla_{\rho}$ is a smooth embedding.

We define the following $2$-form on $Af$ \begin{equation} H(v_1,
v_2)=\int_{U}Tr(v_1\cup v_2)\cup \l^{d-1}, \ \ \ v_1, v_2\in T_a{Af}.
\label{eq:e61} \end{equation} (The integral in question is defined since $\l$ is
a bounded form on $U$, and the $v_j$'s  are $L_2$.) The form $H(.,.)$ is
naturally closed, because it is a constant coefficient $1$-form over an affine
space.

The isomorphism mentioned in Section 4 between $L_2$ and intersection
cohomologies (with coefficients in $\tag$) has the property that corresponding
pairings on degree $1$ cohomology given by (\ref{eq:e32}) and (\ref{eq:e61})
correspond to one another. Thus we get \begin{lem} The pull-back $\beta^*H(.,.)$
is equal to $F$.  \end{lem} The pull-back $\beta^*H(.,.)$ is closed and the
point $\nabla_0$ can be chosen to be an arbitrary representation with trivial
stabilizer, hence $F$ is closed as well.

\

\noindent{\bf Remark.} Let ${\cal A}$ be the space consisting of flat
irreducible connections $\nabla$ on $V$ which are compatible with the unitary
structure and have the property that the difference $\nabla-\nabla_0\in A^1(U,
End({\bar V}))$ is an $L_2$ $1$-form. We also require that $\nabla$ and
$\nabla_0$ have the same monodromies (up to conjugation) around the irreducible
components of $D$.  (It is a mater of simple computations to see that this
condition is redundant and follows from the others.) The tangent space
$T_{\nabla_0}{\cal A}$ is the space of closed unitary $L_2$-forms on $U$ such
that their restriction to a small enough punctured polycylinder
$\Delta^r\times(\Delta^*)^{d-r}$ is $\nabla_0$-exact.  The smooth structure on
${\cal A}$ naturally comes from the fibering $$ \begin{array}{ccc} \G & \to &
{\cal A} \\ {} & {} & \downarrow\pi \\ {} & {} & \MN \end{array}.  $$

It comes down to considering the natural smooth map $s: {\cal A}\to Af$.  The
form $H(.,.)$ on $Af$ is invariant under the action of $\G$ and the pull-back
$s^*H$ of the form $H(.,.)$ to ${\cal A}$ (which we denote by $H_{\cal A}$) is
closed too.

The form $H_{\cal A}$ is $\G$-invariant and vertical, meaning that $H_{\cal
A}(v_1, v_2)=0$, where $v_1, v_2\in T_{\nabla_0}{\cal A}$ and $v_1$ is
$\nabla_0$-exact. This shows that the form $H_{\cal A}$ descends to a $2$-form
on $\MN$ and it can be seen directly from the definitions that this form
coincides with $F$.

\

To finish the proof of Theorem \ref{Th:t11} we need

\begin{lem} The form $F$ is non-degenerate. \end{lem}
\proof  Let $\rho\in\MN$.
The $2$-form $F$ on $T_{\rho}\MN=IH^1(X, \tag)$ is given by
$$ IH^1(X, \tag)\times IH^1(X,
\tag)\to\R:\ \ \ F(x,y)=B\langle x, y\cup\l^{d-1}\rangle. $$
It is clearly skew-symmetric; it is also non-degenerate, because the
Poincar\'{e} duality pairing is non-degenerate and the Hard Lefschetz theorem
says that iterated cap-product by $\l$ induces an isomorphism $$ IH^{d-k}(X
,\tag)\simeq IH^{d+k}(X, \tag).$$  Here we use the fact that the Hard Lefschetz
theorem holds for $L_2$ cohomology \cite{Zu} and we recall (Section 4) the
isomorphism between $L_2$ and intersection cohomologies. $\bigcirc$

\

Combining this result with our above observation that $F$ is closed, we conclude
that $F$ gives a natural symplectic structure on $\MN$.  This proves Theorem
\ref{Th:t11}.

\section{Projective case.}

In this section we assume that $X$ is a projective manifold and we intend to
show that the form $F$ defined in Section 4 is actually a \K form. (One can
conjecture that this is true under the assumption that $X$ is just a \K
manifold, but we do not know the proof of this more general fact.)

  Let $X\subset\C\P^M$ be an embedding of $X$. Let us take a curve $C$ in
$X$ obtained by intersection with $d-1$ generic hyperplanes in $\C\P^M$. This
curve satisfies the following properties. First, $C$ is smooth and, secondly, if
$D=\cup_{i=1}^rD_i$ is the decomposition of $D$ into the union of irreducible
complex analytic subvarieties then $C$ intersects $D_i$ transversally and $C$
does not meet $D_i\cap D_j$ for $i\ne j$. Thus $C\cap D$ is a finite set of
points on $C$. Let $S:=C\min(C\cap D)$ be the open Riemann surface. There is a
natural set of conjugacy classes $\N'$ in the group $G$ assigned to $S$, which
comes from the set $\N$. There is also the moduli space $\MN'$ of flat
irreducible $G$-bundles on $S$ with monodromies around the punctures defined by
the set $\N'$. The following result, which is a direct consequence of Lemma 1.4
in \cite{D3} comes handy. \begin{lem} Let $z\in S$.  The natural morphism
$$\pi_1(S, z)\to \pi_1(U, z)$$ is surjective. \end{lem} This result allows us to
get a smooth inclusion: $$T: \MN \hookrightarrow \MN'.$$

If we pick an irreducible representation $\rho$ of $\pi_1(U)$, it gives rise to
an irreducible representation $\rho'$ of $\pi_1(S)$. There is an obvious
compatibility of images by $\rho$ and $\rho'$ of loops encircling irreducible
components of the divisor $D$ and the punctures on $S$ respectively.

  The tangent space to $\MN'$ is identified (by Proposition \ref{prop:p21})
with the group $Ker[H^1(\pi_1(S), \ag)\to\prod_j H^1({\Gamma_j}', \ag)]$, where
the groups ${\Gamma_j}'\simeq\Z$ are generated by the classes of loops
encircling the punctures.

The above Lemma also gives us the existence of injective homomorphisms:
$$H^1(\pi_1(U), \ag)\hookrightarrow H^1(\pi_1(S), \ag), $$ $$\prod_i
H^1(\Gamma_i, \ag))\hookrightarrow \prod_j H^1({\Gamma_j}', \ag)).$$ (Recall
that $\Gamma_i$ are generated by the classes of loops encircling irreducible
components of the divisor $D$.)  Thus we get the following injection
$$Ker[H^1(\pi_1(U), \ag)\to\prod_i H^1(\Gamma_i, \ag)]\hookrightarrow
Ker[H^1(\pi_1(S), \ag)\to\prod_j H^1({\Gamma_j}', \ag)],$$ which is equal to
$dT:\ T_{\rho}\MN\hookrightarrow T_{\rho'}\MN'.$ (The Lefschetz hyperplane
theorem for intersection cohomology \cite{GMII} gives the same result.) Thus we
proved \begin{prop} The map $T: \MN \to\MN'$ defined above is an embedding.
\end{prop}

Let $F'$ be the $2$-form on the moduli space $\MN'$ defined by Equation
(\ref{eq:e32}) with $d=1$.  Theorem \ref{Th:t11} shows that $F'$ is a symplectic form.
We have the following relation between those two $2$-forms: \begin{prop}
$$F'_{|\MN}=F.$$ \label{prop:p73} \end{prop} \proof Let
$Y\stackrel{\alpha}{\hookrightarrow} X$ be a generic hyperplane section of $X$.
Then we have a well-defined map $$\alpha^*:\ IH^j(X, \tag)\to IH^j(Y, \tag),$$
(which is the transpose to the map $\alpha_*$ of 7.1 \cite{GMII}). Let
$B_X\langle, \rangle$ and $B_Y\langle, \rangle$ be the pairings on $IH^*(X,
\tag)$ and $IH^*(Y, \tag)$ respectively defined as in (\ref{eq:e32}). We have
\begin{lem} Let $\l\in H^2(X)$ be the class of a hyperplane section and let
$a\in IH^j(X, \tag)$, $b\in IH^{d-j-1}(X, \tag)$.  Then $$B_X\langle a,
b\cup\l\rangle =B_H\langle\alpha^*a, \alpha^* b\rangle.$$ \end{lem}

Let $i:C\hookrightarrow X$ be the inclusion and let $[C]=\l^{d-1}\in
H^{2d-2}(X)$ be the class of $C$.  We apply the above result $(d-1)$ times to
  see that if $a\in IH^k(X, \tag)$, $b\in IH^j(X, \tag)$, and $k+j=2$, then
$$B_C\langle i^*a, i^*b\rangle= B_X\langle a,b\cup [C]\rangle.$$ We notice that
when $k=j=1$ the pairings in the left and right hand sides are by definition
$F'$ and $F$. $\bigcirc$

\

  By a well-known theorem of Mehta and Seshadri \cite{MS} the moduli space
$\MN'$ identifies with the moduli space of parabolic stable vector bundles and
hence it has a natural structure of complex manifold. In view of the isomorphism
between intersection and $L_2$-cohomologies, the complex structure on $\MN'$
coming from moduli space of stable parabolic bundles coincides with the almost
complex structure $$J:\  IH^1(C, \tag)\to IH^1(C, \tag)$$ we have introduced in
Section 4.
To conclude that $\MN'$ is a \K manifold, it remains to notice that
if $v,w\in IH^1(X, \tag)$ then $F'(Jv,Jw)=F'(v,w)$.

Moreover, since $C$ is a general curve, the holomorphic $L_2$
$1$-forms on $U$ restrict to holomorphic $L_2$ $1$-forms on $S$ (both with
coefficients in $\tag$). Thus $IH^1(X, \tag)$ is a complex subspace of $IH^1(C,
\tag)$ and Proposition \ref{prop:p73} follows that $\MN$ is a \K submanifold of
$\MN'$.  This ends the proof of Theorem \ref{Th:t12}.

\

\noindent{\bf Remark.} Let us assume that the monodromy around each $D_i$ is of
finite order $k$. Then the same is true for the corresponding monodromy
transformations around the punctures in $S$.  It is well-known
(see e.g. \cite{DW}, \cite{W}) that in this case the form $$\o={kF'\over
4\pi^2}$$ defines an integral cohomology class $[\o]$ in $H^2(\MN', \R)$.
Moreover, there exists a natural holomorphic line bundle $\L'$ on $\MN'$ with
curvature equal to $-2\pi\i [\o]$. Explicit constructions of $\L'$
can be found, for instance, in \cite{DW} and \cite{Kon}.

Using the fact that $\MN$ is a \K submanifold of $\MN'$, we obtain \begin{prop}
Let the monodromy transformation around each $D_i$ be of finite order $k$. Then
there exists a natural line bundle $\L$ on $\MN$ with the curvature equal to $k
F / 2\pi\i$. \end{prop} The line bundle $\L$ is simply the restriction of $\L'$
to $\MN$. The space of $L_2$-sections of $\L$ can serve as a (partial) goal in
the geometric quantization program.

\vskip 0.3in
Department of Mathematics \\ Penn State University \\ University Park, PA 16802
\\ jlb@math.psu.edu, foth@math.psu.edu

\

\noindent{\it AMS subj. class.} \ \ primary 32G13


\thebibliography{123}
\bibitem{Art}{M. Artin, On the solution of analytic equations, {\it Inv. Math.},
{\bf 5}, 1968, 277-291} \bibitem{AB}{M. Atiyah and R. Bott. The Yang-Mills
Equations over a Riemann Surface. {\it Phil. Trans. Roy. Soc}, {\bf A308}, 523
(1982)} \bibitem{A3}{M. Atiyah, {\it The geometry and physics of knots},
Cambridge U.  Press, 1990} \bibitem{AK}{M. Atiyah, {\it K-theory,} Benjamin, New
 York, 1967} \bibitem{BiqTh}{O. Biquard, These, Ecole Polytechnique, 1991}
\bibitem{Biq2}{O. Biquard, Fibr\'es paraboliques stables et connexions
singuli\`eres plates, {\it Bull. Soc. math. France}, {\bf 119}, 1991, 231-257}
  \bibitem{BG}{I. Biswas and K. Guruprasad, Principal bundles on
 open surfaces and invariant functions on Lie groups, {\it Int. J. Math.}, {\bf
 4}, 1993, 535-544} \bibitem{B}{A. Borel {\it at al}. Intersection Cohomology,
Birkh\"{a}user, 1984} \bibitem{BZ}{J.-L.  Brylinski and S. Zucker, An Overview
of Recent Advances in Hodge Theory, {\it Encyclopaedia of Mathematical
Sciences}, {\bf 69}, Several Complex Variables VI, Springer-Verlag, 1990,
39-142} \bibitem{CKS}{E. Cattani, A. Kaplan and W.  Schmid, $L^2$ and
intersection cohomologies for a polarizable variation of Hodge structure, {\it
Invent. Math.}, {\bf 87}, 1987, 217-252}
\bibitem{DW}{G. Daskalopoulos and R. Wentworth, Geometric Quantization for the
Moduli Space of Vector Bundles with Parabolic Structures, preprint, 1992}
\bibitem{De1}{P.  Deligne, Equations Differentielles a Points Singuliers
Reguliers, {\it Lect.  Not. in Math.} {\bf 163}, Springer-Verlag, 1970}
\bibitem{D3}{P. Deligne, Le Groupe Fundamental du
Compl\'{e}ment d'une Courbe Plane n'ayant que des Points Doubles Ordinaires Est
Ab\'{e}lien (d'apr\`{e}s W.  Fulton), S\'{e}m. Bourbaki 543, Nov. 1979, {\it
Lect. Not. Math.}, {\bf 524}} \bibitem{F}{P. A. Foth, Geometry of Moduli Spaces
of Flat Connections on Punctured Surfaces, preprint, {\bf
alg-geom/9703004}, 1996}
\bibitem{GM}{W. Goldman and J. Millson, Deformations of flat bundles over \K
manifolds, {\it Lect. Not. in Pure and Appl. Math.}, {\bf 105}, Dekker, NY,
1987, 129-145} \bibitem{G}{W. Goldman, The Symplectic Nature of Fundamental
Groups of Surfaces, {\it Adv. in Math.}, {\bf 54}, 1984, 200-225}
\bibitem{GMII}{M.  Goresky and R. MacPherson, Intersection Homology II, {\it
Inventiones Mathematicae} {\bf 71}, 1983, 77-129} \bibitem{GHJW}{K. Guruprasad,
J.  Huebschmann, L. Jeffrey and A. Weinstein, Group Systems, Groupoids, and
Moduli Spaces of Parabolic Bundles, preprint, {\bf dg-ga/9510006}, 1995}
\bibitem{KM}{M. Kapovich and J. Millson, The relative deformation theory of
representations and flat connections and deformations of linkages in constant
curvature spaces, {\it Composito Math.}, {\bf 103}, 1996, 287-317}
\bibitem{Kar}{Y.  Karshon, An algebraic proof for the symplectic structure of
moduli space, {\it Proc. AMS}, {\bf 116}, 3, 1992, 591-605} \bibitem{KK}{M.
Kashiwara and T. Kawai, The Poincar\'{e} lemma for variations of polarized Hodge
structures, {\it Publ.  R.I.M.S. Kyoto Univ.}, {\bf 23}, 1987, 345-407}
\bibitem{Kon}{H. Konno, On the Natural Line Bundle on the Moduli Space of Stable
Parabolic Bundles, {\it Comm. Math. Phys.}, {\bf 155}, 1993, 311-324}
\bibitem{MS}{V. Mehta and C. Seshadri, Moduli of Vector Bundles on Curves with
Parabolic Structures, {\it Math. Ann.}, {\bf 248}, 1980, 205-239}
\bibitem{Nits}{N. Nitsure, Moduli of semistable logarithmic connections, {\it J.
Amer. Math. Soc.}, {\bf 6}, 1993, 597-609} \bibitem{Simgar}{C. Simpson, Harmonic
Bundles on Noncompact Curves, {\it J.  Amer. Math. Soc.}, {\bf 3}, 1990,
713-770} \bibitem{SPM}{C.  Simpson, Product of Matrices, {\it Diff. Geom.,
Global Anal. \& Top.}, CMS Conf.  Proc., {\bf 12}, 1992, 157-185} \bibitem{W}{E.
Witten, On Quantum Gauge Theories in Two Dimensions, {\it Comm. Math. Phys.},
{\bf 141}, 1991, 153-209} \bibitem{Zu}{S. Zucker, Hodge theory with degenerating
coefficients:  $L_2$-cohomology in the Poincar\'e metric, {\it Ann. Math.}, {\bf
109}, 1979, 415-476}
\end{document}